 
\input amstex
\documentstyle{gsm}
\NoBlackBoxes
\font\sc=cmcsc10

\font\tit=cmr10 scaled\magstep2
\font\sfa=cmss10 scaled\magstep1
\font\msbm=msbm10
\define\IZ{\hbox{\msbm\char'132}}
 \catcode`\@=11
\def\bqed{\ifhmode\unskip\nobreak\fi\quad
  \ifmmode\blacksquare\else$\m@th\blacksquare$\fi}
\def\mex{\qopname@{mex}}
\catcode`\@=\active
\UseAMSsymbols
\vglue2.1cm
\centerline{\tit Multivision:}
\vskip0.3cm 
\centerline{\tit An Intractable Impartial Game With a 
Linear Winning Strategy}

\vskip1.0cm
\centerline{\sfa Aviezri S.\ Fraenkel}
\vskip0.8cm
\document
{\bf 1. Introduction.} 
Something is definitely wrong. If the game has a linear winning 
strategy, then it is tractable (see Section 3). What's going on? We shall set 
the record straight in just a little while. 

{\sl Multivision\/} is played by two players who move alternately on 
a finite number of piles of finitely many tokens. A move consists of 
selecting a nonempty pile and changing the number of tokens in it. If 
the change is just any {\sl reduction\/} in the number of tokens, we 
have the well-known game of {\sl Nim\/} [1], [2]. In the present case, 
however, it is also permissible, normally, to increase the size of 
a pile by an arbitrarily large factor. 

Specifically, in a multivision game $\Gamma$ we are given $m$ piles 
of tokens, of sizes $n_1,\dots ,n_m$. Let 
$M=M_{\Gamma}=\prod_{i=1}^{m} n_i$. Let $S=\{p_1=2, p_2=3, p_3=5, \dots\}$ 
denote the sequence of primes. If $M>1$, let 
$p_j\in S$ be the smallest prime dividing $M$, and $p_{j+k}$ the 
largest. Thus $p_j<\dots <p_{j+k}$ $(k>0)$, and $\Gamma=\Gamma_{j,k}$
depends on $j$ and $k$. 

A {\sl move\/} consists of selecting any pile $i$ of size greater than 
1, and dividing its size 
$n_i$ by $p_h$ for some prime factor $p_h$ of $n_i$, and multiplying $n_i$ by 
$\prod_{i=h+1}^{j+k} p_{i}^{t_{i}}$ ($\ge 1$), where the $t_i$ 
are arbitrary nonnegative integers. Thus a move consists of reducing 
by 1 the exponent of some prime $p_h$, and increasing the exponents of 
all the primes greater than $p_h$ by arbitrary nonnegative integers. 
Note that if an exponent of a prime $p_h$ $(p_j<p_h\le p_{j+k})$ was 0 
in the initial position or became 0 during play, it may (again) become 
positive, subject to the move rule. 

Play terminates when a position $u$ is reached with all pile sizes 1, so 
no further move is possible. Then $M_u=1$. 
The player reaching this position wins; the opponent loses. Since 
{\it multi\/}plication by a factor at least 1 and di{\it vision\/} by 
a prime is involved in each move, the game is called {\it multivision}. 

Note that we could have defined the game by additions and subtractions 
of tuples of nonnegative integers (exponents) rather than multiplications 
and divisions of products of prime powers. We would then have had to 
define the moves on {\sl vectors\/}, which is less natural than piles. 

\medskip

{\bf Example 1}. Suppose $n_1=\prod_{i=1}^{100} p_i^{p_i}$, 
$n_2=\prod_{i=1}^{100} p_i^{2p_i}$, $n_3=\prod_{i=1}^{100} p_i^{3p_i}$  
($p_i\in S$, $p_{100}=541$). We make a few random initial moves. Player~I 
selects $n_1$, divides $p_{50}^{p_{50}}$ by $p_{50}$ (=229) and replaces  
$n_1$ by say, $n'_1=\prod_{i=1}^{49} p_i^{p_i} p_{50}^{p_{50}-1}
\prod_{i=51}^{100} p_i^{p_i!}$.
Now player~II selects $n_3$, replacing $\prod_{i=1}^{100} p_i^{3p_i}$ by 
say, $\prod_{i=1}^{59} p_i^{3p_i}p_{60}^{3p_{60}-1}
\prod_{i=61}^{100} p_i^{T(3p_i)}$ ($p_{60}=281$), where $T(a)$ is the 
{\it tower function\/}, defined by 
$$T(a)=a^{a^{a^{\cdot^{\cdot^{\cdot^{a^a}}}}}}
\vbox{\hbox{$\Big\}\scriptstyle a$}\kern0pt}.$$
Then player~I selects 
$n'_1=\prod_{i=1}^{100} p_i^{s_i}$, replacing it by say, 
$\prod_{i=1}^{36} p_i^{s_i}p_{37}^{s_{37}-1}\prod_{i=38}^{100}
p_i^{T^{(p_i)}(p_i)}$, 
and so on, where $T^{(a)}$ is the $a$-th iterate of the tower function 
$T$.\medskip

This example shows that play can continue for quite a while before it 
terminates---if it does terminate!\medskip

The questions we are interested in:
\roster
\item Does every play of every multivision game terminate, or are there 
plays that continue forever?\par
\item What's a winning strategy for the game if and when it exists?\par
\item What's the complexity of such a winning strategy? That is, how 
difficult is it to: (i) discover a winning strategy, (ii) compute a 
winning move, and (iii) consummate a winning strategy, i.e., actually 
reach the end of play by winning the game?\par
\endroster

We begin by establishing a small framework.\medskip 

{\bf Definitions and Notations.} 
\roster
\item"$\bullet$" A {\it follower\/} (or {\it option\/}) of a game position 
$u$ is any position that is reachable from $u$ in a single move.\par 
\item"$\bullet$" The set of all followers of position $u$ is denoted by 
$F(u)$.\par
\item"$\bullet$" A game position $u$ is an $N$-{\it position\/} if 
the {\it N}ext player can force a win, i.e., the player moving from $u$.\par 
\item"$\bullet$" A game position $u$ is a $P$-{\it position\/} if 
the {\it P}revious player can force a win, i.e., the opponent of the 
player moving from $u$.\par
\item"$\bullet$" The set of all $N$-positions of a game is denoted by 
$\Cal N$ and the set of all $P$-positions by $\Cal P$.
\endroster
\medskip

Evidently, the goal of each player is always to move from an $N$-position 
into a $P$-position.\medskip 

The main purpose of this note is to prove

\proclaim{\bf Theorem 1} Every play of every multivision game\/ $\Gamma$ 
terminates. Play at any position\/ $u$ of\/ $\Gamma$ can be prolonged 
arbitrarily long precisely as long as\/ $M_u$ has a prime factor 
$\ne p_{j+k}$, where\/ $M_u=\prod_{i=j}^{j+k} p_i^{s_i}$.\endproclaim

In other words, every play of $\Gamma$ terminates, but it can be 
prolonged arbitrarily long precisely while $M_{\Gamma}$ is divisible by 
a prime less than $p_{j+k}$. 

\proclaim{\bf Theorem 2} Every multivision game has a winning strategy 
for precisely one of the two players. The\/ $P$-positions are precisely the 
set of positions\/ $u$ for which\/ $M_u$ is a square.\endproclaim 

Note that the initial position of the game in Example~1 is a 
$P$-position.\medskip

Summarizing, the 2-player game multivision has the following properties:
\roster
\item"(i)" precisely one of the two players can consummate a win in 
a finite number of moves,\par
\item"(ii)" at each stage the winner can compute a next winning move 
in linear time,\par
\item"(iii)" at each step, except the last few, play can be 
delayed indefinitely,\quad and\par
\item"(iv)" one can compute in linear time who of the two players can 
win and who loses, assuming the winner plays correctly (which can easily 
be done by (ii)). 
\endroster 

So the conclusion is that the game does have a definite winner, that is, a 
player who can force a win in a finite number of moves. Moreover, the winner's 
winning moves can be computed linearly and we can determine the winner 
in linear time. Yet the game is intractable, since (iii) implies that 
the winner would have to live indefinitely to consummate a win. Thus 
the intractability part of the title is the truth, and the linear 
strategy is a lie, since only {\sl parts} of the strategy are linear. 
A strategy is tractable only if all its parts are (see Section 3). 

A 2-player game $\Gamma$ is {\it impartial\/} if, for every position of 
$\Gamma$, both players have the same options. Otherwise $\Gamma$ is 
{\it partizan\/} [1]. Thus Nim and multivision are impartial, and 
chess and Go are partizan. Impartial games are usually simpler than 
partizan games. Yet even the former may exhibit pathological behavior, 
such as demonstrated here.\bigskip  
 
{\bf 2. Proofs.}\nolinebreak 
\demo{\bf Proof of Theorem 1} With any position $u$ of the game associate 
a $(k+1)$-tuple of nonnegative integers $(t_j,\dots ,t_{j+k})$, which 
are the exponents of $M_u=\prod_{i=j}^{j+k} p_i^{t_i}$, $t_i\ge 0$; 
also say that $M_u$ or $(t_j,\dots ,t_{j+k})$ {\it encodes\/} 
the position $u$. Thus the set of all game positions corresponds to a 
subset $K$ of all $(k+1)$-tuples of nonnegative integers. 

Consider the lexicographic ordering $\prec$ of $K$:
$$(r_j,\dots ,r_{j+k})\prec (t_j,\dots ,t_{j+k})$$ 
if there exists $h\in\{j,\dots ,j+k\}$ such that $r_i=t_i$ for all $i<h$, and $r_h<t_h$. 
Suppose that $(t_j,\dots ,t_{j+k})$ encodes 
position $u$. Any follower $v$ of $u$ is obtained by reducing by 1 the 
exponent of some $p_h$ in some $n_i$, and possibly increasing the 
exponents of some of the primes $>p_h$. Therefore $v$ is encoded by 
$(t_j,\dots t_{h-1},t_h-1,r_{h+1},\dots ,r_{j+k})$, where $r_i\ge t_i$ for 
$i\ge h+1$. Thus $v\prec u$. Since the lexicographic ordering is a 
well-ordering of the set $K$, play terminates after a finite 
number of moves. 

If $u$ is a position encoded by $M_u=p_{j+k}^{s_{j+k}}$, then for $r>s_{j+k}$, 
play terminates in fewer than $r$ moves. 
But if $M_u=ap_i^{s_i}p_{j+k}^{s_{j+k}}$ where $p_i<p_{j+k}$,\ \ $s_i>0$ 
and all the prime factors of $a$\ ($a\ge 1$) are less than $p_i$, 
then there exists a pile 
of size $a'p_i^{s'_i}p_{j+k}^{s'_{j+k}}$ with $s_i'>0$ and $a'$ a factor 
of $a$. The move to $a'p_i^{s'_i-1}p_{j+k}^{s'_{j+k}+r}$ 
realizes play of length $\ge r$ moves.\bqed\enddemo\medskip

\demo{\bf Proof of Theorem 2} It follows from Theorem~1 of [4] 
that for any game which terminates in a finite number of moves, the 
set of all its positions can be partitioned uniquely into subsets 
$\Cal N$ and $\Cal P$. From Theorem~3 of [4], it follows that for 
{\sl acyclic\/} games, such as multivision, any partition of 
the game positions into $\Cal N'$ and $\Cal P'$ for which 
$$u\in\Cal N'\quad \text{if and only if}\quad F(u)\cap\Cal P'\ne\emptyset\quad 
\text{and}\quad u\in\Cal P'\quad \text{if and only if}\quad 
F(u)\subseteq\Cal N',$$
satisfies $\Cal N'=\Cal N$ and $\Cal P'=\Cal P$. 
It thus suffices to show: 
\roster
\item"(i)" Any position $u$ encoded by $M_u=\prod p_i^{s_i}$ such that $M_u$ 
is not a square, {\sl has\/} a follower $v$ encoded by a square.\par
\item"(ii)" {\sl Every\/} follower $v$ of a position $u$, encoded by 
$M_u=\prod p_i^{s_i}$, which is a square, is encoded by a nonsquare. 
\endroster\medskip

(i)\ If $M_u$ is not square, it has a smallest prime factor, say $p_h$, 
with an odd exponent $s_h$. Then there exists a pile of size $n_i$, such 
that $p_h^s$ divides $n_i$, $p_h^{s+1}$ does not divide $n_i$ and $s$ 
is odd. The move 
$n_i\rightarrow n_ip_h^{-1}\prod_{i=h+1}^{k} p_i^{t_i}$ results 
in a position $v$ encoded by $M_v$ such that $M_v$ is a square, if and 
only if $t_i$ has the same parity as the exponent of $p_i$ in $M_u$. 

(ii)\ Any move $n_i\rightarrow n_ip_h^{-1}\prod_{i=h+1}^{k} p_i^{t_i}$, 
from a position encoded by a square, 
results in a position $v$ encoded by a nonsquare $M_v$, since the exponent of 
$p_h$ in $M_v$ is odd.\bqed\enddemo\medskip

Theorem~1 states, roughly, that play can be prolonged arbitrarily 
long except during the ``end moves''. Delaying the end of play is 
normally in the interest of the loser. Nevertheless, we show explicitly 
that either player can effect a delay.\medskip 

\proclaim{\bf Corollary} Play at any position\/ $u$ of a multivision 
game\/ $\Gamma$ can be prolonged indefinitely by either player, 
precisely while it is true that\/ $M_u$ has a prime factor\/ $<p_{j+k}$, 
where\/ $M_u=\prod_{i=j}^{j+k} p_i^{s_i}$.\endproclaim

\demo{\bf Proof} We first show that the winner can prolong play 
indefinitely. Let $M_u=ap_i^{s_i}p_{j+k}^{s_{j+k}}$, where 
$p_i<p_{j+k}$,\ \ $s_i>0$ and all the prime factors of $a$ are 
less than $p_i$. Suppose first that $M_u$ is square. By 
Theorem~2, player~II can force a win. In the worst case player~I 
keeps reducing $s_{j+k}$. After it becomes 0, player~I is the first to reduce 
$s_i$ or the exponent of a prime less than $p_i$. Then player~II responds by 
making the exponent of $p_{j+k}$ arbitrarily large, while at the same time 
restoring $M$ to be square. If $M_u$ is not a square, then 
player~I can win by moving $u\rightarrow v$ with $M_v$ a square, so 
we are back in the previous case. The same argument shows that the loser 
can delay play indefinitely.\bqed\enddemo\bigskip

{\bf 3. Complexity of Multivision.} 
It is convenient to input the initial set of pile sizes $n_i$ in 
the form of their prime decompositions. During play it is convenient 
to maintain for each $n_i$ and for $M$ a $(k+1)$-vector $V_i$ and
$V_M$ over $\IZ_2$, where $V_M$ is the sum over $\IZ_2$ of the $V_i$. 
If $n_i=\prod_{r=j}^{j+k} p_r^{s_r}$, then $V_i$ is the vector 
$(s_j,\dots ,s_{j+k})$ over $\IZ_2$. 
Thus for Example~1, the vectors have initially length 100, 
$V_1$ and $V_3$ consist of 1s except for the leftmost 
component, which is 0 in both, and $V_2$ and $V_M$ are the 0-vector. Whenever 
the loser makes a move, $V_M$ becomes nonzero.  The winner then 
locates the leftmost column of $V_M$, say $h$, which contains a 
1. There exists some $n_i$ such that $V_i$ contains a 1 in column 
$h$, so the winner can restore $V_M$ to 0 by dividing $n_i$ by 
$p_h$ and multiplying it by $p_r$ raised to an odd (even) power for 
all $r>h$ for which $V_M$ has a 1 (0) in column $r$.

In addition to the binary vectors $V_i$, also their true values 
(over the ring of integers) have to be maintained for checking the 
validity of the move of the opponent and for recognizing the end of 
the game. 

If $f$ and $g$ are functions from the nonnegative integers into the 
nonnegative reals, we write, as customary, $f(n)=\Theta (g(n))$ if 
there exist constants $c_2\ge c_1>0$ such that $c_1g(n)\le f(n)\le c_2g(n)$ 
for all large $n$. For complexity studies it's usually preferable to use 
the $\Theta$-notation rather than the $O$-notation, since if, say, 
$f(n)=O(\log n)$, then also $f(n)=O(n)$, so the $O$-notation doesn't 
necessarily differentiate between functions one of which is exponential 
in the other. 

In the realm of discrete mathematics we say, roughly, that a problem 
is {\it tractable\/} if it can be solved in time (number of steps) 
that is a polynomial in the input length of the problem. Otherwise it is 
{\it intractable}. A game is tractable if: (i) the recognition 
of the $P$-positions, (ii) the computation of a winning move, and 
(iii) the consummation of the win are tractable. It is 
reasonable to define tractability of (i) and (ii) as for all other 
problems in discrete mathematics, namely computation in time polynomial 
in the input size. But it is unreasonable to do so for (iii): 
Nim is the prototype of a tractable game. Given a game of Nim 
with two piles, each of length $n$. The input length is then 
$\Theta (\log n)$, since we can represent $n$ in binary, say, which 
requires $\lceil\log_2 (n+1)\rceil$ bits (binary digits). If player~I 
keeps taking a single token from one of the 
piles, player~II has to take one from the other pile in order to maintain 
a winning position, so the length of play is $\Theta (n)$, which is 
exponential in the input length. But a game with play length beyond 
a polynomial in an exponential is intractable. See [3, Section 7] for 
more information about the tractability of games. 

The winning strategy for multivision described in the proof of Theorem~2 
is linear, i.e., tractable with polynomial of degree 1, and is similar to 
the winning strategy of Nim. Theorem~2 and 
its proof show that, in fact, (i) and (ii) of question~(3) can be done 
in linear time. Yet the game is intractable, because the loser can 
normally prevent the winner's consummation of the well-deserved win 
indefinitely. Example~1 perhaps begins to suggest that play 
can last beyond anything suggested by the Ackermann function, which 
is the prototype of an extremely fast-growing function; see, e.g., 
[5, Section 5], [10, Section 2]. 

Summarizing this section, we may say that multivision is a simple 
example of an intractable game---which has a linear winning strategy
if we disregard the length of play. There are games for which (i) of 
question~(3) is provably hard (exponential), and others for which (i) 
and (iii) are easy but (ii) is either undecidable or its complexity 
is unknown. Multivision is a simpler game than the long Epidemiography 
games [5], [6], [7], which are related to the Hercules-Hydra game (reviewed 
in [11]), yet it can last much longer.\medskip 

We end by describing a generalization $K$-{\sl multivision\/} of multivision. 
Let $K$ be a fixed integer $\ge 2$. We begin to play $K$-multivision 
with $m$ piles of sizes $n_1,\dots ,n_m$ as in multivision. A 
{\sl move\/} consists of selecting a prime $p_h$ and an integer $s<K$ 
and dividing some piles $n_{i_1},\dots ,n_{i_r}$ by 
$p_h^{s_1},\dots ,p_h^{s_r}$ with $\sum_{i=1}^{r} s_i=s$  (where $r\le s$; 
and $n_{i_1}$ is divisible by $p_h^{s_1},\dots ,n_{i_r}$ by $p_h^{s_r}$; 
$s_i$ positive integers). We also multiply $n_i$ by 
$\prod_{i=h+1}^{k} p_{i}^{t_{i}}$ ($\ge 1$), where the $t_i$ are 
arbitrary nonnegative integers, for $i\in \{i_1,\dots ,i_r\}$.  

The special case $K=2$ of $K$-multivision is, of course, multivision. 

\proclaim{\bf Theorem 3}  Every play of every $K$-multivision game\/ $\Gamma$ 
terminates. Play at any position\/ $u$ of\/ $\Gamma$ can be prolonged 
arbitrarily long, even by the winner, precisely as long as\/ 
$M_u$ has a factor $<p_{j+k}$, where $M_u=\prod_{i=j}^{j+k} p_i^{s_i}$. 
Every $K$-multivision game has a winning strategy for precisely one of 
the two players. The\/ $P$-positions are precisely the 
set of positions\/ $u$ for which\/ $M_u$ is a $K$-th power.\endproclaim 

The proof is similar to that of Theorem~1; it follows from the following 
two easily verifiable observations:

\roster
\item"(i)" Any position $u$ encoded by $M_u=\prod p_i^{s_i}$ such that $M_u$ 
is not a $K$-th power, {\sl has\/} a follower $v$ encoded by a 
$K$-th power.\par
\item"(ii)" {\sl Every\/} follower $v$ of a position $u$, encoded by 
$M_u=\prod p_i^{s_i}$, such that $M_u$ is a $K$-th power, is encoded by 
an integer $M_v$ such that $M_v$ is not a $K$-th power. 
\endroster
Finally we remark that other pathological games are discussed, e.g., 
in [8], [9], and [13].\medskip

{\bf Confession.} While reading the section  on factoring a number $x$ 
via the difference of squares method, $x=y^2-z^2$, in the inspiring 
survey article [12] of Carl Pomerance, multivision formed before my eyes 
in a flash. The sad part is that I then promised myself to finish 
reading [12], which is still in the bottomless queue of things I promised 
myself to do. Instead I wrote this paper that same day and sent it, on 
an impulse, to some friends, which I immediately regretted. But some of 
my recipients, among them Thomas Ferguson, Daniel Ullman, and Herbert Wilf, 
sent me nice feedback, and Herb encouraged me to send it to the {\sc Monthly}. 
So if you don't like this paper, blame me. If you do like it, thank Tom, 
Dan, and Herb.\bigskip 

\centerline{\bf References}\medskip

1. E.\ R.\ Berlekamp, J.\ H.\ Conway, and R.\ K.\ Guy, {\it Winning
Ways\/} (two volumes), Academic Press, London, 1982.

2. C.\ L.\ Bouton, Nim, a game with a complete mathematical theory,
{\it Ann.\ of Math.} {\bf 3} (1902) 35--39. 

3. A.\ S.\ Fraenkel, Scenic trails ascending from sea-level Nim to 
alpine chess, 
in: {\it Games of No Chance\/}, Proc. MSRI Workshop on Combinatorial 
Games, July, 1994, Berkeley, CA (R. J. Nowakowski, ed.), MSRI Publ. Vol.~29, 
Cambridge University Press, Cambridge, 1996, pp. 13--42.

4. A.\ S.\ Fraenkel, Combinatorial game theory foundations applied to 
digraph kernels, {\it Electron. J.\ Combin.\/} {\bf 4} (2) 
(1997) \#R10.\quad 
http://www.combinatorics.org/Volume\_4/wilftoc.html

5. A.\ S.\ Fraenkel, M.\ Loebl, and J.\ Ne\v set\v ril, 
Epidemiography---II. Games with a dozing yet winning
player, {\it J.\ Combin. Theory\/} (Ser.\ A) {\bf 49} (1988) 129--144.

6. A.\ S.\ Fraenkel and M.\ Lorberbom, Epidemiography with various 
growth functions, {\it Discrete Appl.\ Math.\/} {\bf 25} (1989) 53--71. 

7. A.\ S.\ Fraenkel and J.\ Ne\v set\v ril, Epidemiography, 
{\it Pacific J.\ Math.\/} {\bf 118} (1985) 369--381. 

8. J.\ P. Jones, Some undecidable determined games, {\it Internat.
J. Game Theory\/} {\bf 11} (1982) 63--70.

9. J.\ P. Jones and A.\ S. Fraenkel, Complexities of winning
strategies in diophantine games, {\it J. Complexity\/} {\bf 11} (1995) 435-455.

10. J.\ Ne\v set\v ril, Ramsey Theory, Ch.~25 in: {\it Handbook of
Combinatorics\/} (R.\ L.\ Graham, M.\ Gr\"otschel, and L.\ Lov\'asz, eds.), 
Elsevier, Amsterdam, 1995. 

11. J.\ Ne\v set\v ril and R.\ Thomas, Well quasi ordering, long games 
and combinatorial study of undecidability, {\it Contemp.\ Math.\ 
Proc.\ Symp. AMS\/} {\bf 65} (1987) 281--293. 

12. C.\ Pomerance, A tale of two sieves, {\it Notices AMS\/} 
{\bf 43} (1996) 1473--1485.    

13. M.\ O.\ Rabin, Effective computability of winning strategies,
Contributions to the Theory of Games vol.~3, {\it Ann. of Math. Stud.} 
{\bf 39}, 147--157, Princeton, 1957.\medskip 

Department of Applied 
Mathematics and Computer Science, Weizmann Institute of Science, 
Rehovot 76100, Israel\medskip

{\tt fraenkel\@wisdom.weizmann.ac.il}

\enddocument